\documentclass{amsart}
\usepackage[all]{xy}
\usepackage{graphicx}
\usepackage{psfrag}
\usepackage{amsthm}
\usepackage{amscd}
\usepackage{amsfonts}
\usepackage{amstext}
\usepackage{amssymb}
\usepackage{amsmath}
\usepackage{amsxtra}
\usepackage{plain}
\usepackage[all]{xy}

\newcommand{\comm}[1]{}

\numberwithin{equation}{section}
\let\newpf\proof  
\newenvironment{pf}{\newpf\proofname}{\qed\endtrivlist}

\let\newpf\proof

\def\be{\begin{equation}}
\def\ee{\end{equation}}

\def\ba{\begin{align}}
\def\ea{\end{align}}

\newtheorem{prop}{{\bf Proposition}}

\newtheorem{thm}[prop]{{\bf Theorem}}

\newtheorem{rem}{\bf Remark}

\newcommand{\Z}{\mathbb{Z} }
\newcommand{\N}{\mathbb{N} }

\newcommand{\R}{\mathbb{R} }

\newcommand{\ind}{\rm{ind}}

\newcommand{\A}{\mathcal{A}}

\newcommand{\Ch}{\rm{Ch}}

\newcommand{\End}{\rm{End}}

\newcommand{\STr}{\rm{STr }}
\newcommand{\cl}{\ensuremath{\operatorname{cl}}}

\title{Heat equation approach to index theorems on odd dimensional manifolds}
\author{Mostafa Esfahani Zadeh
\address{Institute for Advanced Studies in 
Basic Sciences(IASBS), Zanjan-IRAN, {\bf and}, Mathematisches Institut 
Georg-August-Universit\"at G\"ottingen-Germany}}
\email{esfahani@iasbs.ac.ir}
\begin{document}
\maketitle

\begin{abstract}
 D.Freed has formulated and proved an index theorem on odd dimensional 
spin manifolds with boundary. 
The proof is based on analysis by Calderon and Seeley. 
In this note we are going to give a proof of this theorem 
using the heat kernels methods for boundary
conditions of Dirichlet and Von Neumann type. Moreover we consider also the 
Atiyah-Patodi-Singer 
spectral boundary condition which is not considered in Freed's paper.
As a direct consequence of the method, we will obtain some 
information about isospectral invariants of the boundary conditions.
This proof does not uses the cobordism invariance of index and are easily
generalized to family case. 
\end{abstract}
\section{Introduction}Dirac operators on compact odd dimensional 
manifolds, without boundary, are formally self adjoint, so their index will be 
zero and have no interest. However, if the underlying manifold has 
boundary, by puttin suitable boundary conditions, one can obtain 
non self adjoint operators with non zero index. In \cite{Freed1}, 
D. Freed, inspired by the work of physicians has formulated a such 
index theorem and proved it by means of symbol calculus of elliptic 
boundary problems. In this paper we give a heat equation proof of this 
theorem and consider, in addition, the case of spectral boundary conditions. 
To compute the contribution of APS condition, we use the perfect symmetry 
between positive and negative part of the spectrum of Dirac operators in 
even dimension rather than using calculus which are 
established to handle this problem in general (see, eg. \cite{MePi} 
and \cite{BiCh1}).
In section \ref{secdo} we state the theorem and prove it for one 
Dirac operator. At the end of this section we will give a necessary 
condition for isospectrality of these boundary problems 
with exchanged boundary conditions on different connected components 
of boundary. In section \ref{sec3} we will formulate and prove the 
index theorem for a family of Dirac operators by studying the Chern character of 
a superconnection adopted to the family of Dirac operators. 
Although our notation refers to the Dirac operator acting on standard 
complex spinor fields, the method can be applied directly to the case 
of Dirac type operators acting on twisted Dirac operator or Dirac 
operators acting on Clifford bundles. 
\newline 
{\bf Acknowledgment}:The author would like to thank Pierre B\'erard 
and Thomas Schick for useful discussions and suggestions.
\section{Index theorem for a Dirac type operator}\label{secdo}
 Let $M$ be a compact spin manifold of dimension $2n+1$ with a 
Riemannian metric  taking the product form 
 $d^2u+g$ in a collar neighborhood $U=[0,1^+)\times\partial M$ of the 
 boundary. Here $1^+$ may be any real number greater
 than $1$. In the sequel we denote the coordinate of the interval $[0,1^+)$ 
by $u$ and $v$ and those of $\partial M$ by $y$ and $z$. 
A typical point of $M$ will be denoted by $x$. We fix a spin structure 
on $M$ giving rise to the complex spin 
vector bundle $S$ on $M$. 
The orientation and spin structure on $M$ induce an orientation 
and a spin structure on $\partial M$. So the
 restriction of spin bundle $S$ to the collar neighborhood of the 
boundary splits into the positive and negative parts
\[S_{|U}\simeq S^+\oplus S^-\]
where $S^+\oplus S^-$ is the spin bundle of $\partial M$ lifted to the 
collar neighborhood $U$ in the obvious manner.
According to this splitting, each spinor field $\phi\in C^\infty(M,S)$ has a 
decomposition $\phi_{|U}=\phi^+\oplus\phi^-$, where
 $\phi^\pm\in C^\infty(\partial M,S^\pm)$. With respect to this 
decomposition, 
the Dirac operator $D$ has the following form where $A$  denote the Dirac 
operator of $\partial M$ (see, eg. \cite[chapter 9]{BoWo1}) 
\[
\left(\begin{array}{cc}i\partial_u&iA^-\\-iA^+&-i\partial_u\end{array}\right).
\]
Let $\partial M=\sqcup_iN_i$ where $N_i$'s are disjoint union of connected 
components of boundary. For each $i$, choose $\epsilon_i$
be $0$, $+$ or $-$ arbitrarily and fix them. For $\epsilon_i=\pm$ let 
$P^{\epsilon_i}$  denote the following local boundary condition
\begin{equation}\label{loccon}
(\phi_{|N_i})^{\epsilon_i}=0.
\end{equation}
Corresponding to $\epsilon_i=0$, $P^0$ denotes the Atiyah-Patodi-Singer(APS) 
boundary condition. We recall briefly its definition. 
The Dirac operator $A_i$ of each connected component $N_i$ of $\partial M$ 
has a discrete resolution $\{\phi_\lambda,\lambda\}$. 
Since $N_i$ is even dimensional, there is a symmetry 
betweed positive and negative part of spectrum. This symmetry is given by 
the following unitary isomorphism $\mathcal U$ between eigenspaces with 
opposite eigenvalues 
\begin{equation}\label{sympn}
\mathcal U(\phi^+_\lambda+\phi_\lambda^-)=\phi^+_\lambda-\phi_\lambda^- ~~;~~
\phi_\lambda=\phi^+_\lambda+\phi_\lambda^-.
\end{equation}
 This symmetry is responsible for the eta-invariant 
of $N$ vanishes. Let $P^0\in\End(L^2(N_i.S))$ be the orthogonal projection on 
the subspace generated by $\phi_\lambda$'s 
with $\lambda\geq0$. This projection defines the APS boundary condition. 
By $P^\epsilon$ we denote the boundary condition which equals 
$P^{\epsilon_i}$ on $N_i$.  
\newline
Let $\phi$ and $\psi$ be smooth spinor field on $M$, then the following 
formula holds
\begin{equation}\label{adrel}
<D\,\phi,\psi>_{L^2}-<\phi,D\,\psi>_{L^2}=-\int_{\partial 
 M}<\cl(\frac{\partial}{\partial u})\phi,\psi>\,dy.
\end{equation}
So it turns out that boundary problems $(D,P^+)$ and 
$(D,P^-)$ are adjoint, one for other. Moreover the adjoint problem 
for $(D,P^0)$ is $(D,Id-P^0)$. 
We denote the formal adjoint of the boundary problem 
$D:=(D,P^\epsilon)$ by $D^*:=(D,P^{\bar\epsilon})$. 
These boundary problems are elliptic so they have finite 
dimensional kernels consisting of smooth spinor fields. 
Therefore one can define the index of this problem by 
\begin{equation}\label{ind1}
\ind(D,P^\epsilon)=\dim\ker\,D^*D-\dim\ker\,DD^*.
\end{equation}
\begin{thm}[See theorem B of \cite{Freed1}]\label{thm1}
Denote by $A_i$ the Dirac operator of the even dimensional manifold $N_i$. 
The following formula holds 
\[\ind(D,P^\epsilon)=\frac{1}{2}\sum_{\epsilon_i=-}\ind\,A_i-\frac{1}{2
}\sum_{\epsilon_i=+}\ind\,A_i-\frac{1}{2}\sum_{\epsilon_i=0}\dim\ker A_i.\]
\end{thm}
To prove this theorem, we consider these boundary problems on  
half cylinder $\R^{\geq0}\times N$ with product spin structure and 
product Riemannain metric. Here $N$ may denote each one of $N_i$'s and $A$ 
denotes its Dirac operator. 

\subsection{Index density of local boundary conditions on 
half cylinder}\label{avalin}
We consider the boundary condition $P^+$; the case $P^-$ 
can be studied in the same manner. 
On this half cylinder the Dirac operator $D$, 
acting on compactly supported spinor fields subjected to condition $P^+$, 
takes the following form 
\begin{equation}\label{dircyl}
\left(\begin{array}{cc}i\partial_u&iA^-\\-iA^+&-i\partial_u\end{array}\right).
\end{equation}
We denote the adjoint problem by $D^*$.  
The Laplacian operators $D^*D$ and $DD^*$, 
acting on compactly supported spinor fields $\phi(u,y)$,
 are non negative operators taking the form $-\partial^2_u+A^2$. 
The induced boundary condition for $D^*D$ is 
\[\phi^+(0,y)=0~~\text{ and }~~(\partial_u\phi^-+A^+\phi^+)_{|u=0}=0.\]
Since $A$ is a tangential operator, these conditions reduce to Dirichlet 
condition for $\phi^+$ and Von Neumann 
condition for $\phi^-$
\begin{gather}
\phi^+(0,y)=0\label{Dirboucon}\\
\frac{\partial \phi^-}{\partial u}(0,y)=0.\label{Vonboucon}
\end{gather}
Concerning the operator $DD^*$, a similar argument shows that the 
boundary conditions take the following forms
\begin{equation}\label{adboucon}
\frac{\partial}{\partial u}\psi^+(0,y)=0~~\text{ and }~~\psi^-(0,y)=0.
\end{equation}
Consider, for $t>0$, the heat operators $e^{-tD^*D}$. 
The kernel of this operator, $\bar K_1(t,u,v,y,z)$, with respect 
to the boundary 
condition \eqref{Dirboucon} has the following explicite form, 
cf. \cite{McSi1}
\begin{gather}\label{con11}
\frac{1}{\sqrt{4\pi t}}\lbrace 
\exp(\frac{-(u-v)^2}{4t})-\exp(\frac{-(u+v)^2}{4t})\rbrace 
e^{-tA^-A^+}(t,y,z),
\end{gather}
while the heat kernel for boundary condition \eqref{Vonboucon} is 
\begin{equation}\label{con2}
\frac{1}{\sqrt{4\pi t}}\lbrace 
\exp(\frac{-(u-v)^2}{4t})+\exp(\frac{-(u+v)^2}{4t})\rbrace 
e^{-tA^+A^-}(t,y,z).
\end{equation}
So the trace density of heat operator $e^{-tD^*D}$, as a function of 
$t>0$ and $u\geq0$ is 
\begin{align}
\bar K_1(t,u):&=\int_N\rm{tr_y}\bar K_1(t,u,v,y,y)\,dy\notag\\
&=\frac{\rm{Tr}\,e^{-tA^-A^+}}{\sqrt{4\pi t}}\lbrace 
1-\exp(\frac{-u^2}{t})\rbrace
+\frac{\rm{Tr}\,e^{-tA^+A^-}}{\sqrt{4\pi t}}
\lbrace 1+\exp(\frac{-u^2}{t})\rbrace.\label{hetden}
\end{align}
Denote by $\bar K_2(t,u,v,y,z)$ the kernel of heat operator $e^{-tDD^*}$ 
with boundary conditions \eqref{adboucon}. Above discussion gives the 
following expression for trace density of this heat operator  
\begin{equation}\label{adhetden}
\bar K_2(t,u)=\frac{\rm{Tr}\,e^{-tA^-A^+}}{\sqrt{4\pi t}}\lbrace 
1+\exp(\frac{-u^2}{t})\rbrace
+\frac{\rm{Tr}\,e^{-tA^+A^-}}{\sqrt{4\pi t}}
\lbrace 1-\exp(\frac{-u^2}{t})\rbrace.
\end{equation}
So we obtain the following formula for trace density
$\bar K_+(t,u)=\bar K_1(t,u)-\bar K_2(t,u)$ of the operator 
$e^{-tD^*D}-e^{-tDD^*}$ with boundary condition $P^+$
\begin{equation}\notag
\bar K_+(t,u)=\frac{e^{-\frac{u^2}{t}}}{\sqrt{\pi t}}\lbrace 
Tr\,e^{-tA^+A^-}-Tr\,e^{-tA^-A^+}\rbrace.
\end{equation}
Integrating with respect to $u\in\R^{\geq0}$ we get 
\begin{equation}\label{filcyl}
\int_0^\infty \bar K_+(t,u)\,du=-\frac{1}{2}\ind A.
\end{equation}
Here we have used the McKean-Singer formula \(\ind 
A=Tr\,e^{-tA^-A^+}-Tr\,e^{-tA^+A^-}\).
Similarly, if $\bar K_-(t,u)$ denote the trace density of operator 
$e^{-tD^*D}-e^{-tDD^*}$ with boundary condition $P^-$, then 
\begin{equation*}
\int_0^\infty \bar K_-(t,u)\,du=\frac{1}{2}\ind A.
\end{equation*}
If instead of integration on $[0,\infty)$ in above we integrate  
on the finite interval $[0,\frac{1}{2}]$, the difference, 
would be of exponential decay when $t$ goes toward $0$.
In fact we have the following relation 
\begin{equation}\label{adfilcyl}
\int_0^{\frac{1}{2}}\bar K_\pm(t,u)\,du=\mp\frac{1}{2}
\ind A+O(e^{-\frac{1}{4t}})
\end{equation}
Following relation is a result of above discussion and we write down it 
here for future reference 
\begin{equation}\label{abed1}
\bar K_\pm(t,u)\sim0~~\text{ exponentially at }t=0 \text{ for }u\neq0.
\end{equation}
\begin{rem}\label{clrem}
If $y\neq z$, then it is well known that the heat kernel $e^{-tA^-A^+}(t,y,z)$ 
and $e^{-tA^+A^-}(t,y,z)$, as well as thir derivatives
 with respect to $y$, are exponentially small at $t=0$. If $y=z$ and $u\neq v$ 
the expression given in formulas 
\eqref{con11} and \eqref{con2} have this property. In this case the 
differentiation may be taken with respect to $t$.
\end{rem}

\subsection{Index density of APS boundary problem on half 
cylinder}\label{ssloc} 
Now let $\epsilon_i=0$ and consider the half cylinder $\R^{\geq0}\times N$ 
with boundary problems $D=(D,P^0)$ and 
$D^*=(D,Id-P^0)$. The induced boundary condition for $D^*D$ is 
\[P^0(\phi)=0~~\text{ and }~~(Id-P^0)D\phi=0.\]
Exchanging the role of $P$ and $(Id-P)$ we obtain the adjoint induced boundary 
condition $Q^*$ for $D^*D$.
Denote by $E_\lambda$ the product vector bundle on $\R^{\geq0}$ whose fibers 
are the $\lambda$-eigenspace of $A$. Operators $D^*D$ and $DD^*$ 
takes the form $-\partial^2_u+\lambda^2$ on sections of $E_\lambda$.  
Put $A_\lambda:=A_{|E_\lambda}$. 
Boundary conditions on section $\phi(u,y)$ of $E_\lambda$, for $D^*D$, read 
\begin{align}
&\phi(0,y)=0&\text{ for }\lambda\geq0~;\label{aps1}\\
&(\frac{\partial}{\partial u}+\lambda)_{|u=0}\,\phi(u,y)=0 
&\text{ for }\lambda<0.\label{aps2}
\end{align}
Following boundary conditions must be considered for $DD^*$.  
\begin{align}
&\psi(0,y)=0&\text{ for }\lambda<0~;\label{aaps1}\\
&(\frac{\partial}{\partial u}+\lambda)_{|u=0}\,\psi(u,y)=0 
&\text{ for }\lambda\geq0.\label{aaps2}
\end{align}
Let $K_1^{\lambda}$ and $K_2^{\lambda}$ denote, respectively, the  
heat kernels of $D^*D$ and $DD^*$. 
Let $\bar K(t,u)$ denote the supertrace density 
\footnote{Here Tr and STr denote 
the $L^2$ trace while tr and str denote the pointwise trace of finite dimensions vector 
 space}
$\rm{Tr}_N\,e^{-tD^*D}(t,u)-\rm{Tr}_N\,e^{-tDD^*}(t,u)$. 
We are interested in the following quantity which is in fact the index of $D$, 
acting on spinor fields supported in $[0,1/2)\times N$ and subjected to 
APS condition at $\{0\}\times N$ 
\begin{equation}\label{eigsum}
\int_0^{\frac{1}{2}}\bar K(t,u)
=\sum_\lambda\int_0^{\frac{1}{2}}K_1^\lambda(t,u)- 
\int_0^{\frac{1}{2}}K_2^{-\lambda}(t,u).
\end{equation}
We recall operator  $\mathcal U$ provides, for $\lambda\neq0$, a 
unitary isomorphism 
\[\mathcal U_{\lambda}:=\mathcal U:C^\infty(\R^{\geq0},E_\lambda)\to 
C^\infty(\R^{\geq0},E_{-\lambda})\] 
One has $-A_{-\lambda}=\mathcal U_\lambda A_\lambda\mathcal U_\lambda^{-1}$, so 
$DD^*_{|E_{-{\lambda}}}
=\mathcal U_\lambda D^*D_{|E_{\lambda}}\mathcal U_{\lambda}^{-1}$. 
Moreover $\mathcal U_\lambda$ exchanges the boundary conditions \eqref{aps1} 
and \eqref{aaps1}. So in expression \eqref{eigsum}, the terms indexed by 
$\lambda>0$ cancel themselves. 
The trace density of fundamental solution of 
$\partial_t-\partial^2_u+\lambda^2=0$ with 
boundary condition \eqref{aps2}  is given by 
following expression, cf. \cite[relation 2.17]{APS1}
\begin{equation*}\label{hetap1}
K_1^\lambda(t,u)=\left(\frac{e^{-\lambda^2t}}{\sqrt{4\pi t}}\lbrace 
1+\exp(\frac{-u^2}{4t})\rbrace 
+\lambda e^{-2\lambda u}
\rm{erfc\lbrace\frac{u}{\sqrt{t}}-\lambda\sqrt{t}\rbrace}\right)
.\dim E_\lambda, ~~\lambda<0. 
\end{equation*}
The fundamental solution with respect to boundary condition \eqref{aaps2} is 
\begin{equation*}\label{hetap2}
K_2^{-\lambda}(t,u)=\left(\frac{e^{-\lambda^2t}}{\sqrt{4\pi t}}\lbrace 
1+\exp(\frac{-u^2}{4t})\rbrace 
-\lambda e^{2\lambda u}
\rm{erfc\lbrace\frac{u}{\sqrt{t}}+\lambda\sqrt{t}\rbrace}\right)
.\dim E_{-\lambda}, ~~\lambda<0. 
\end{equation*}
The error function is defined by the following formula 
\[\rm{erfc}(x)=\frac{2}{\sqrt{\pi}}\int_x^\infty e^{-s^2}\,ds.\]
Subtracting above expressions, we get the following 
relations at $t=0$ for $\lambda<0$ 
\begin{align*}
K_1^\lambda(t,u)-K_2^{-\lambda}(t,u)&
\sim 0~~~\text{ exponentially for }u\neq0,\\ 
K_1^\lambda(t,0)-K_2^{-\lambda}(t,0)&\sim\lambda+o(1).
\end{align*}
So 
\begin{equation*}
\int_0^{\frac{1}{2}}K_1^\lambda(t,u)-K_2^{-\lambda}(t,u)\,du\sim 0~~
\text{ exponentially at }t=0.
\end{equation*}
Therefore, when $t$ goes zero, the only terms 
with probably nonzero contribution  in relation \eqref{eigsum} are 
indexed by $\lambda=0$.
For $\lambda=0$, the boundary conditions \eqref{aps1} and \eqref{aaps2} 
are respectively the Dirichlet and Von Neumann boundary conditions. 
Regarding expression \eqref{con11} and \eqref{con2} we have 
\begin{equation}\label{laze}
K_1^0(t,u)-K_2^0(t,u)=-\frac{e^{\frac{-u^2}{t}}}{\sqrt{\pi t}}\dim\ker A.
\end{equation}
Hence 
\[\int_0^{\frac{1}{2}}(K^1_0(t,u)-K^2_0(t,u))\,du
\sim -\frac{1}{2}\dim\ker A~~\text{ exponentially at }t=0.\]
We summarize above discussion in following proposition
\begin{prop}\label{ahah}
Denote the heat kernel of boundary problems $D^*D$ and $DD^*$, respectively, 
by $\bar K_1$ and $\bar K_2$.
 With above notation, put $E^+:=\oplus_{\lambda>0}E_\lambda$ and 
$E^-:=\oplus_{\lambda<0}E_\lambda$. $\mathcal U$ provides a  
natural unitary isomorphism between these spaces and satisfies 
\begin{align*}
\bar K_{2|E^-}&=\mathcal U \bar K_{1|E^+}\mathcal U^{-1}\\
\bar K_{2|E^+}&\sim\mathcal U \bar K_{1|E^-}\mathcal U^{-1}
+o(1)~~\text{ exponentially at }t=0. 
\end{align*} 
Moreover, with $\bar K(t,u):=\bar K_1(t,u)-\bar K_2(t,u)$ we have 
\begin{equation}\label{apscyl}
\int_0^{\frac{1}{2}}\bar K(t,u)\,du
\sim-\frac{1}{2}\dim\ker A~~\text{ exponentially at }t=0.
\end{equation}
\end{prop}
\begin{rem}
For later use we rewrite the following result of above discussion
\begin{equation}\label{abed}
\bar K(t,u)\sim 0~~\text{ exponentially for } u\neq0.
\end{equation}
\end{rem}

\subsection{Index theorem for one Dirac operator }\label{prforone}  
The heat equation proof of the theorem \ref{thm1} is based on the 
Mckean-Singer  
formula. Let $D$  denote the boundary problem $(D,P^{\epsilon})$ and 
$D^*$ denote 
its adjoint problem $(D,P^{\bar\epsilon})$. These boundary 
problems are elliptic, so the elements of $\ker(D,P^\epsilon)$ are smooth 
spinor fields which are exactly those elements of $\ker D^*D$  
satisfying  induced boundary conditions. Similar remark applies to 
$\ker(D^*,P^{\bar\epsilon})$ and $\ker DD^*$. 
Hence using the spectral
resolution of $D^*D$ and $DD^*$ one has the following Mckean-Singer 
type equality 
\[\ind(D,P^\epsilon)=\rm{Tr}\,e^{-tD^*D}-\rm{Tr}\,e^{-tDD^*}~;~~t>0.\]
Denote by $\bar K_{1\epsilon}$ and $\bar K_{2\bar \epsilon}$, respectively, 
the fundamental solutions of the cylindrical heat operators 
$e^{-tD^*D}$ and $e^{-tDD^*}$.
Let $K_\epsilon$ and $K_{\bar\epsilon}$ be, respectively, the fundamental 
solution of $e^{-tD^*D}$ and $e^{-tDD^*}$. We 
are going to give asymptotic expressions, at $t=0$, for these fundamental 
solutions in term of  $\bar K_{1\epsilon}$, $\bar K_{2\bar\epsilon}$ and of 
the fundamental solution $\tilde K$ of $e^{-tDD}$ on the 
double of $M$, i.e. $M\sqcup_{\partial M}M^-$.
For this purpose, following \cite[Page 54]{APS1}, let $\rho(a,b)$ be a smooth 
increasing function on $\R^{\geq0}$ 
such that 
\[\rho(u)=0 \text{ for }u\leq a~~;~~\rho(u)=1\text{ for }u\geq b.\]
Collar neighborhoods of connected component of boundary are assumed being 
parameterized by $u\in[0,1^+]$, so the following 
functions may be considered as smooth function on $M$ with constant extension into $M$.
\begin{align*}
&f_2=\rho(\frac{1}{4},\frac{1}{2}),~~~~~~~&g_2=\rho(\frac{1}{2},\frac{3}{4})\\
&f_1=1-\rho(\frac{3}{4},1),~~~~~~~&g_1=1-g_2
\end{align*}
Put 
\begin{equation}\label{initst}
\mathcal K_\epsilon=f_1\bar K_{1\epsilon} g_1+f_2\tilde K g_2~~\text{ and }~~\mathcal K_{\bar\epsilon}=f_1\bar K_{2\bar\epsilon} g_1+f_2\tilde K g_2.
\end{equation}
Since $f_i=1$ on the support of $g_i$, one conclude that $\mathcal K_\epsilon$, 
as an operator on $C^\infty(M,S)$, goes toward $Id$ when $t\to0$.
Moreover the remark \ref{clrem} show that $(\frac{\partial}{\partial t}+D^*D) 
\mathcal K_\epsilon$ is exponentially small, 
out of diagonal, when $t\to0$. These two condition are sufficient for 
using $\mathcal K_\epsilon$ as the initial step in 
construction of heat kernel using the Levi's sum, cf. \cite{McSi1}. As a 
consequence, the difference between heat kernel
$K_\epsilon$ and $\mathcal K_\epsilon$ is exponentially small when $t$ 
goes toward $0$. This argument applies also to 
$K_{\bar\epsilon}$ and $\mathcal K_{\bar\epsilon}$. Therefore we have 
\begin{align*}
\ind(D,P^\epsilon)&=\int_{diag(M)}\lbrace tr\,K_\epsilon(t,x,x)-tr\, 
K_{\bar\epsilon}(t,x,x)\rbrace\\
&=\lim_{t\to0}\int_{diag(M)}\lbrace tr\,\mathcal K_\epsilon(t,x,x)-tr\, 
\mathcal K_{\bar\epsilon}(t,x,x)\rbrace\\
&=\lim_{t\to0}\sum_i\int_0^{\frac{1}{2}}\bar K_{1\epsilon}(t,u)-
\bar K_{2\bar\epsilon}(t,u)\\
&=\lim_{t\to0}\sum_i\int_0^{\frac{1}{2}}\bar K_{\epsilon}(t,u)\,du.
\end{align*}
To deduce the last equality, we have used the fact that the contribution 
of the trace of the heat
operator on double $M\sqcup M^-$ is the same in the expressions 
$tr\,\mathcal K_\epsilon$ and
$tr\,\mathcal K_{\bar\epsilon}$. We have use also $f_1(u)=g_1(u)=1$ for $0\leq 
u\leq\frac{1}{2}$. 
Now relations \eqref{filcyl}, \eqref{adfilcyl} and 
 \eqref{apscyl} imply the desired formula 
\[\ind(D,P^\epsilon)=\frac{1}{2}\sum_{\epsilon=-}
\ind A_i-\frac{1}{2}\sum_{\epsilon=+}\ind A_i-\frac{1}{2}
\dim\ker (A_i).\]

\subsection{Isospectrality problem}
In this subsection we denote by $Q^\epsilon$ the second degree boundary 
condition associated to $P^\epsilon$.
As we have mentioned before, for $t>0$, the heat operator 
$e^{-tD^*D}$ is a compact 
smoothing self adjoint operator, so the boundary condition problem 
$(D^*D,Q^\epsilon)$ 
has a resolution $\{(\phi_\lambda,\lambda)\}$. Let $\epsilon'$ 
be an another set of boundary conditions. For simplicity let $\epsilon_i$ 
and $\epsilon'_i$  be $+$ or $-$. Let 
\begin{gather*}
e^{-tA_i^-A_i^+}\sim t^{-n}\sum_{k=0}^\infty 
a^+_{ik}t^k~,\\
e^{-tA_i^+A_i^-}\sim t^{-n}\sum_{k=0}^\infty a^-_{ik}t^k~
\end{gather*}
be the asymptotic expansions of heat operators on $N_i$. It follows from 
Mckean-Singer formula that 
\begin{equation}\label{mcsi}
a^+_{ik}=a^-_{ik}~~\text{ for }k\neq n~\text{ and }
a^+_{in}=a^-_{in}+\ind A_i.
\end{equation}
If the boundary problems 
$(D^*D, Q^\epsilon)$ and $(D^*D,Q^{\epsilon'})$ have the same spectrum 
 then for $t>0$ the trace difference  
$Tr(e^{-tD^*D},Q^\epsilon)-Tr(e^{-tD^*D},Q^{\epsilon'})$ vanishes. 
We are interested to study the asymptotic behavior, at $t=0$, of this 
trace difference to get some spectral invariant of these kind of boundary 
problems. 
Relations \eqref{initst} and discussion following them show that this 
trace difference is asymptotic to the following expression at $t=0$ 
\[\int_M\mathcal K_\epsilon(t,x,x)-
\int_M\mathcal K_{\epsilon'}(t,x,x)=\int_0^{1/2} 
\int_{\partial M}\bar K_{1\epsilon}(t,y,u)
-\bar K_{1\epsilon'}(t,y,u)\,dy\,du.\]
On the other hand, relations 
\eqref{hetden} and \eqref{adhetden} can be used to deduce that the last 
expression 
is in its turn asymptotic to 
\begin{equation*}
(\frac{1}{(4\pi t)})^{\frac{n+1}{2}}\sum_{k,i}\lbrace \sum_{\epsilon_i=+} 
t^ka^-_{ik}+\sum_{\epsilon_i=-} t^ka^+_{ik}-\sum_{\epsilon'_i=+} 
t^ka^-_{ik}-\sum_{\epsilon'_i=-}t^ka^+_{ik}\rbrace. 
\end{equation*} 
In view of relations \eqref{mcsi} this last expression simplifies to the 
following one
\[\sum_{\epsilon_i=-,\epsilon'_i=+}\ind 
A_i-\sum_{\epsilon_i=+,\epsilon'_i=-}\ind A_i.\] 
So, using the cobordism invariance of index $\sum_{\epsilon_=\pm}\ind A_i=0$ 
, we get the following theorem
\begin{thm}\label{isos}
With notation of theorem \ref{thm1}, if two boundary problem 
$(D^*D,Q^\epsilon)$ and $(D^*D,Q^{\epsilon'})$ are isospectral then 
\[\sum_{\epsilon_i=-,\epsilon'_i=+}\ind 
A_i=\sum_{\epsilon_i=+,\epsilon'_i=-}\ind A_i=0.\]
Moreover this condition is the only one which can be deduced from heat 
equation asymptotic formulas.
\end{thm}
\begin{rem}
Comparing this result with the case of scalar Laplacian should be 
interesting.
Let $(M,\partial M)$ be a smooth manifold with $\partial M=N_0\sqcup N_1$. 
Let $\triangle$ be the scalar Laplacian acting on smooth 
functions which are subjected to Dirichlet condition on $N_0$ and Von Neumann 
condition on $N_1$. It follows directly from relations \eqref{con11} 
and \eqref{con2}
 that all heat equation invariants\footnote{This means the invariant 
coming from asymptotic expansion, at $t=0$, of $Tr\,e^{-t\triangle}$} of 
$A_{|N_1}$ are spectral invariants for $\triangle$ with this mixed boundary 
conditions. In particular $Vol(N_1)$, so $Vol(N_0)$, are included in the 
spectral invariants. Thought this may seem in contrast with above theorem, 
notice that our boundary conditions are at the same time of Dirichlet 
and Von Neumann type.
\end{rem}
\begin{rem}\label{prob1}
Let $P^+$ denote $\epsilon_i=+$ for all $i$ and define $P^-$ similarly.
Using the cobordism invariance of index, 
the condition of above theorem is satisfied for these pure boundary 
conditions. This leads to the following interesting question:  
Does there exist a spin manifold with boundary such that $(D^*D,Q^+)$ 
and $(D^*D,Q^-)$ be isospectrum.
\end{rem}

\section{Index theorem for families}\label{sec3}
Theorem \ref{thm1} can be generalized to include families of Dirac type
operators. At first we recall geometric setting for family index theorems.
\newline 
Let $M\hookrightarrow F\to B$ be a fibration of odd dimensional compact spin 
manifolds, with 
boundary, over a compact smooth manfold $B$. The boundaries of fibers form an 
another fibration $\partial M\hookrightarrow F'\to B$. 
This boundary fibration has a fibred collar 
neighborhood $U$ of the form $[0,1^+)\times F'\subset F$ 
which restricts to a collar neighborhood in each fiber. We denote by 
$\pi$ the projection on the second factor. Assume that the fibration $F$ 
is endowed with a smooth fiberwise 
Riemannian metric which is of product form $d^2u+g$ in the 
collar neighborhood $U$. 
Here $g$ is a fiberwise metric on the boundary fibration 
$F'$. We assume
also that the fiberwise spin structure is of product form in $U$. 
By these assumption, a typical
fiber $M$ of the fibration $F$ satisfies conditions described in 
previous section. So in the sequel
we will use freely the previous notation in this family context.
Let $D$ and $A$ denote, respectively, the fiberwise family of Dirac 
operators associated to fibrations $F$ and $ F'$. 
Put $ F'=\sqcup_iN_i$ where each $N_i$ is a connected fibration over 
$B$ whose fibers are
even dimensional compact spin manifolds. For each $i$, let $\epsilon _i$ be 
$0$, $+$ or $-$ arbitrarily and fix it. The family 
$(D,P^\epsilon)$ of boundary problems determines 
an analytic index $[\ind(D,P^\epsilon)]$  in
$K_0(B)$. In other hand, the boundary family $A$ determines, in its turn, 
a class $[\ind A]$ in $K_0(B)$(see \cite{AtSi4}).
The equality of these classes is announced in \cite{Freed1}. 
Here we prove the following theorem using the heat equation tools applied 
to superconnections.
\begin{thm}\label{thm2}
The following equality holds in $H^*_{dr}(B)$
\[\Ch[\ind(D,P^\epsilon)]=\frac{1}{2}\sum_{\epsilon_i=
-}\Ch[\ind A_i]-\frac{1}{2}\sum_{\epsilon_i=+}\Ch[\ind A_i]
-\frac{1}{2}\sum_{\epsilon_i=0}\Ch(\ker A_i).\]
\end{thm}
For being able to use the language of the theory of super graded 
differential 
modules, we consider the direct sum $S\oplus S$ as a 
smooth family of super spin bundle such that the first and the second 
summand are respectively the even and 
the odd part. The Clifford action of a vertical tangent vector $v\in TM_b$ 
on $S_b\oplus S_b$ is given by the following grading reversing matrix 
\[\left(\begin{array}{cc}
0&\cl(v)\\ \cl(v)&0
\end{array}\right).
\] 
With this Clifford action, we obtain the following families of grading 
reversing Dirac operators on $F$ and on $F'$
\[\mathsf D:=\left(\begin{array}{cc}
0&D^*\\D&0
\end{array}
 \right)~;~~\mathsf A:=\left(\begin{array}{cc}
0&A\\A&0
\end{array}
 \right)
\]
Here $D$ denotes the boundary problem $(D,P^\epsilon)$ while $D^*$ denote 
$(D,P^{\bar\epsilon)}$.
These operators act on vertical spinor fields $\phi\oplus\psi$ 
satisfying  boundary 
conditions $P^\epsilon(\phi)=0$ and $P^{\bar\epsilon}(\psi)=0$.
For each $b\in B$, $\mathsf D_b$ is a vertical self adjoint odd differential 
operator. The induced second degree boundary problem is $\mathsf D^2$ with 
boundary conditions \eqref{Dirboucon}, \eqref{Vonboucon} and \eqref{adboucon}.
In below we deal with the infinite dimensional bundle $\mathcal E$ of vertical
spinor fields. A typical fiber $\mathcal E_b$ of this bundle is
$C^\infty(M_b,S)$. Bundle $\mathcal E\oplus \mathcal E$ is $\Z_2$-grading in 
obvious way. 
To prove theorem \ref{thm2} we need  
a connection $\nabla$ on sections of $\Lambda^*B\otimes\mathcal E$. 
\begin{rem}\label{rem3}
In following two subsections, we will study the case of a fibration 
with half cylindrical fibers endowed with some 
more special connections. Using a partition of unity on half line $\R^{\geq0}$ 
we assume that the connection $\nabla$ coincides with theses special 
connections in collar neighborhood $[0,1^+)\times F'$. 
\end{rem} 
We need to extend the action of $D$ and $A$ on sections 
$\omega\otimes\xi$ of $\Lambda^*B\otimes\mathcal E$. 
The extension are given by following relations
\[D(\omega\otimes\xi)=(-1)^{\deg\omega}\omega\otimes D\xi~;~
A(\omega\otimes\xi)=(-1)^{\deg\omega}\omega\otimes A\xi.\]
Now $\nabla\oplus\nabla$, $\mathsf D$ and $\mathsf A$ 
 are graded differential operators acting on smooth sections of  
graded bundle $\Lambda^*B\otimes\mathcal E\oplus\Lambda^*B\otimes\mathcal E$. 
In bellow we shall study the superconnections  
$\mathbb B=\mathsf D+\nabla\oplus\nabla$ 
adapted to the family $\mathsf D$ of Dirac operators and the connection 
$\mathbb A=A+\nabla$ which is adopted to family $A$. 
We denote by $\mathbb B_t$ and $\mathbb A_t$ their rescaled versions. 
The rescaled curvature 
$\mathbb F_t:=\mathbb B_t^2$ has the form $t\mathsf D^2+\mathbb F_{t[+]}$ 
where $\mathbb F_{t[+]}$ is the following  differential operator with 
differential form coefficients of positive degree 
\begin{equation}\label{curv}
\mathbb F_{t[+]}=:\left(\begin{array}{cc}
\nabla^2&t^{1/2}(D\nabla+\nabla D)\\t^{1/2}(D\nabla+\nabla D)&\nabla^2
\end{array}\right).
\end{equation}
We shall be dealing with the structure of heat operator of rescaled curvature 
, so we explain briefly how to construct it.
Let $e^{-t\mathsf D^2}=e^{-tD^*D}\oplus e^{-tDD^*}$ be the heat operator of 
$(\mathsf D,P^\epsilon\oplus P^{\bar\epsilon})$
and let $R$ be a family of smoothing operator. 
The heat kernel of perturbed family $\mathsf D_s:=\mathsf D+sR$, 
for $0\leq s\leq1$, is given by Voltera formula: 
\begin{gather*}
e^{-t\mathsf D_s}=e^{-t\mathsf D^2}+\sum_{k=1}^{\infty}(-t)^k
I_k(t\mathsf D_s,sR),\\
I_k(t\mathsf D,sR):=\int_{\triangle_k}e^{-s_0t\mathsf D^2}
\,sR\,e^{-s_1t\mathsf D}\,sR\dots e^{-s_{k-1}t\mathsf D^2}
\,sR\,e^{-s_kt\mathsf D^2},\\
\triangle_k=\{(s_0,s_1,\dots,s_k)\in\R^{k+1}|s_i\geq0~;~\sum_is_i=
1.\}\end{gather*}
Because $R$ is smoothing, the operator $e^{-t\mathsf D}R$ has a smooth kernel 
for $t\geq0$ and $\|e^{-t\mathsf D}R\|_\ell\leq C(\ell)\|R\|_\ell$. So 
\[\|I_k(t\mathsf D,R)\|_\ell\leq\frac{C(\ell)^{k+1}\|R\|_\ell^k}{k!},\]
which implies the convergence of above sum in $C^\ell$-norm. But it 
is clear from above Voltera formula that 
\begin{equation}\label{nimeh}
e^{-t\mathsf D_s}-e^{-t\mathsf D}=o(1),\quad\text{ at }t=0,
\end{equation}
so the smoothing perturbation 
$sR$ has no effect on asymptotic behavior of heat operator at $t=0$. 
Let $\mathbb B_s:=\mathbb B+sR$ be the superconnection adopted to 
perturbed family $\mathsf D_s$. Its rescaled curvature has the form 
$\mathbb F_{s,t}=t\mathsf D_s+\mathbb F_{s,t[+]}$, where 
\begin{equation}\label{dileh}
\mathbb F_{s,t[+]}=\mathbb F_{t[+]}+O(t).
\end{equation}
Heat operator of supercurvature $\mathbb F_{s,t}$ is given again by  
Voltera formula 
\begin{gather}
e^{-\mathbb F_{s,t}}=e^{-t\mathsf D_s^2}+\sum_{k=1}^{\dim B}(-1)^k
I_k(t\mathsf Ds,\mathbb F_{s,t[+]}),\label{voltera}
\end{gather}
It should be clear from this construction that $e^{-\mathbb F_{s,t}}$ 
is a vertical family  of smoothing operator with coefficients 
in $\Lambda^*(B)$. So its supertrace is finite and 
defines an element in $\Omega^*(B)$. Although  
$I_k(t\mathsf D_s,\mathbb F_{s,t[+]})$ 
depends on involved operator in a rather complicated way, its asymptotic 
behavior at $t=0$ is simple to describe. 
At first it follows from \eqref{nimeh} and \eqref{dileh} that 
$I_k(t\mathsf D_s,\mathbb F_{s,t[+]})
=I_k(t\mathsf D,\mathbb F_{t[+]})+o(1)$.
Moreover, in expression $I_k(t\mathsf D,\mathbb F_{t[+]})$ 
the contribution of off-diagonal operators in \eqref{curv} can be neglected, 
at $t=0$, because these operators are multiplied by $t^{1/2}$. 
So $I_k(\mathsf D,\mathbb F_{t[+]})
=I_k(\mathsf D,\nabla^2\oplus\nabla^2)+o(1)$ and we get the following 
relation that we will use later. This relation 
holds in $\Omega^*(B)$ with respect to $C^\ell$-norms
\begin{equation*}
e^{-\mathbb F_{s,t}}= e^{-t\mathsf D^2+\nabla^2\oplus\nabla^2}+o(1).
\end{equation*}
The Chern character of the rescalled perturbed superconnection is defined by 
$\Ch(\mathbb B_{s,t}):=\rm{STr}\,e^{-\mathbb F_{s,t}}$. 
To investigate the relation between $\Ch(\mathbb B_{s,t})$ and 
$\Ch(\mathbb B_{t})$ we recall a general theorem about superconnections 
an super bundles, cf. \cite[Theorem 9.17]{BeGeVe}. Let $\mathsf B_\tau$ be a 
smooth family of superconnection on a differential superbundle.  
The following formula holds 
\cite[Theorem 9.17]{BeGeVe})
\begin{equation}\label{rasta}
\frac{d}{d\tau}\Ch(\mathsf B_\tau)
=-d\,\rm{STr}\,(\frac{d\,\mathsf B_\tau}{d\tau}
e^{-\mathsf B_\tau^2})\in\Omega^*(B).
\end{equation}
If we apply this formula to $s$-dependent family of superconnections 
$\mathbb B_{s,t}$ we get 
\[\Ch(\mathbb B_{1,t})-\Ch(\mathbb B_{t})
=-d\int_0^1\STr(\frac{d\,\mathbb B_{s,t}}{ds}\,e^{-\mathbb F_{s,t}})\,ds\]
Therefore the perturbation $R$ does not affect the class of Chern character 
in de-Rham cohomology of $B$. 
Impact of this perturbation on behavior at $t=\infty$ may 
be crucial. In fact there is a general methods (see \cite[Lemma 2.1]{AtSi4}) 
to construct a self adjoint perturbation $R$ such that 
$\dim\ker(\mathsf D+R)$ be independent of $b\in B$. 
In this case $\ker(\mathsf D+R)$ is a smooth 
finite dimensional super vector bundle over $B$, so it determines a
class in $K^0(B)$. This class, being independent of the perturbation, 
is denoted by $[\ind\,D]$ and is called the analytical index of family 
$D$. We summarize these discussion in following proposition
\begin{prop} We can assume that $\dim\ker(\mathsf D_b)$ is independent 
of $b\in B$ by perturbing by self adjoint smoothing operators. 
This perturbation does not affect the class of Chern form $\Ch(\mathbb B_t)$ in 
$H^*_{dr}(B)$ not its behavior when $t$ goes toward $0$. Moreover one has 
 \begin{equation}\label{asvel}
e^{-\mathbb F_{t}}= e^{-t\mathsf D^2+\nabla^2\oplus\nabla^2}+o(1).
\end{equation}
With this assumption we have $[\ind\,D]=[\ker\,\mathsf D]\in K^0(B)$. 
\end{prop}
Our proof of theorem \ref{thm2} is based on a precise study of the 
behavior of $\Ch(\mathbb B_t)$ at $t=0$ and $t=\infty$ and comparing them.
Following above proposition, $\ker\mathsf D$ is 
a vector bundle on $B$ and the formal difference of its even and odd parts i.e, 
$\ker(D,P^\epsilon)-\ker(D^*,P^{\bar\epsilon})$ represents the index  
class $[\ind\,D]\in K^0(B)$. Let $Q_0$ be the projection on 
$\ker \mathsf D$ which is a smooth
family of vertical smoothing operators. It is clear that 
$\nabla_0=Q_0\nabla Q_0$ is a connection 
on the vector bundle $\ker \mathsf D$. Therefore the differential form 
$\rm{Str}\,e^{-\nabla_0^2}$ is closed 
and provides a representation of $\Ch(\ind\,D)\in H^*_{Dr}(B)$. 
\begin{prop}\label{pro1}
\begin{enumerate}
\item The following convergence occurs in $\Omega^*(M)$ with respect to 
uniform $C^\ell$-norm for each $\ell\in\N$
\[\lim_{t\to+\infty}\Ch(\mathbb B_t)=\Ch(\ker\,\mathsf D,\nabla_0).\]
\item For $t>0$, the Chern form $\Ch(\mathbb B_t)$ is closed and its class 
in de-Rham cohomology $H^*_{dr}(B)$ is independent of $t$, in particular 
\[\Ch(\mathbb B_t)=\Ch(\nabla_0)\in H^*_{dr}(B)~~\text{ for }~t>0.\]
\end{enumerate}
\begin{pf}
Consider the following orthogonal decomposition of superbundles 
\begin{equation}\label{delstan}
\mathcal E\oplus\mathcal E=\ker\,\mathsf D\oplus\rm{Im}\,\mathsf D.
\end{equation}
Heat operator $e^{-t\mathsf D}$ is a family of smoothing 
non negative operators parameterized by the compact set $B$. So there is a 
uniform gap around $0$, in the spectrum of each element of this family. 
This simple observation and general properties of graded 
nilpotent algebras can be used to get the following relation with respect 
to above direct sum decomposition. 
(see \cite[page 290]{BeGeVe} ) 
\begin{equation}\label{hakem}
e^{-\mathbb B_t}\sim\left(\begin{array}{cc}
e^{-\nabla_0^2}&0\\0&0
\end{array}
 \right)
+\left(\begin{array}{cc}
O(t^{-1/2})&O(t^{-1/2})\\O(t^{-1/2})&O(t^{-1})
\end{array}
 \right).
\end{equation}
Notice that the convergence with respect to $C^\ell$-norm is implicit in 
above asymptotic formula which proves the first part of the proposition.  
To prove the second part we apply relation \eqref{rasta} to family 
$\mathbb B_t$. So for $t_2,t_1>0$  
\begin{equation*}
\Ch(\mathbb B_{t_2})-\Ch(\mathbb B_{t_1})
=-d\int_{t_1}^{t_2}\rm{STr}\,(\frac{d\,\mathbb B_t}{dt}
e^{-\mathbb F_t})\,dt.
\end{equation*}
If we regard this relation in de-Rham cohomology group $H^*_{dr}(B)$ we 
obtain the desired relation in the second part of the proposition.
\end{pf}
\end{prop}
\begin{rem}
If $t_1$ goes to $\infty$ in above relation then one can deduce the 
following stranger relation
\begin{equation}\label{delst}
\Ch(\mathbb B_{t_2})-\Ch(\nabla_0)
=d\int_{t_2}^{\infty}\rm{STr}\,(\frac{d\,\mathbb B_t}{dt}
e^{-\mathbb F_t})\,dt\in\Omega^*(B),
\end{equation}
provided the integral in the right hand side is finite. To prove 
the convergence of this integral, we notice that  
\[\frac{d\,\mathbb B_t}{dt}
=\frac{1}{2t^{1/2}}\left(\begin{array}{cc}
0&0\\\mathsf D&0
\end{array}
 \right),\]
with respect to decomposition \eqref{delstan}. This relation with \eqref{hakem} 
give rise to 
\[\rm{STr}\,(\frac{d\,\mathbb B_t}{dt}
e^{-\mathbb F_t})=O(t^{-3/2}).\]
So the above integral is convergent at $t=\infty$.
\end{rem}
Now we shall study the behavior of $\Ch(\mathbb B_t)$ 
when $t\to0$. For this purpose we give, as in the previous section,  
an another description of heat operator
$\,e^{-\mathbb F_t}$ in term of heat operator on double of $F$ and heat 
operator on half cylinder fibration $\R^{\geq0}\times F'$ with typical 
fiber $\R^{\geq0}\times\partial M$. As before, let $\pi$ denote the projection 
on second summand.
All local structures in this cylindrical case are exactly the same of the 
collar neighborhood $U$, so, for example, the Dirac operator takes the 
form \eqref{dircyl}. 

\subsection{Index density of local conditions on family of half 
cylinders} 
Let $\mathcal E_0=\mathcal E_0^+\oplus\mathcal E_0^-$ 
be the bundle over $B$ whit typical fiber 
$C^\infty(\partial M_b,S^+\oplus S^-)$. 
Let $\nabla_0=\nabla_0^+\oplus\nabla_0^-$ be a connection on $\mathcal E_0$. 
Then $\bar\nabla=\pi^*\nabla_0$ is a connection on 
$\mathcal E:=\pi^*\mathcal E_0$ which preserves the grading. 
We denote by $\bar{\mathbb B}$ the superconnection 
$\mathsf D+\bar\nabla\oplus\bar\nabla$ acting on $\mathcal E\oplus\mathcal E$. 
We consider $\R^{\geq0}\times N$ 
with local boundary condition $P^+$. Here $N$ may denote each one of 
$N_i$ and for simplicity we drop the index $i$. From relation \eqref{asvel}; 
as far as we are interested in the asymptotic behavior of heat operators 
at $t=0$; we can replace $\bar{\mathbb F}_t:=\bar{\mathbb B}_t^2$ by 
fallowing operator 
\begin{equation}\label{curcy2}
\left(\begin{array}{cc}
-t\partial^2_u&0\\0&-t\partial^2_u
\end{array}
\right)+\pi^*\left(\begin{array}{cc}
tA^2+\nabla^2_0&0\\0&tA^2+\nabla^2_0
\end{array}
\right)
\end{equation}
Relation \eqref{asvel} can be applied to superconnection 
$\mathbb A=A+\nabla_0$ where, $\nabla_0=\pi^*\nabla^+_0\oplus\pi^*\nabla^-_0$. 
This provides $e^{-\mathbb A_t^2}= e^{-(tA^2+\nabla^2_0)}+o(1)$, so 
\begin{align}
e^{-\bar{\mathbb F}_t}= e^{-\bar{\mathsf  F}_t}+o(1)\label{asa},
\end{align}
where 
\begin{equation*}
\bar{\mathsf  F}_t^2:=\left(\begin{array}{cc}
-t\partial^2_u&0\\0&-t\partial^2_u
\end{array}
\right)+\pi^*\left(\begin{array}{cc}
\mathbb A_t^2&0\\0&\mathbb A_t^2
\end{array}
\right).
\end{equation*}
Notice that the operator $\partial_u$ commutes with all
other operators involved in above expression, so the results of 
previous section can be used to give explicite expression for  
$e^{-\bar{\mathsf  F}_t}$. For example its 
even part, acting on $\phi^+(u,y)\oplus\phi^-(u,y))$ subjected to 
boundary conditions \eqref{Dirboucon} and \eqref{Vonboucon} is 
\[\frac{1}{\sqrt{4\pi t}}\left(
\begin{array}{cc}
\exp(\frac{-(u-v)^2}{4t})-\exp(\frac{-(u+v)^2}{4t})&0\\
0&\exp(\frac{-(u-v)^2}{4t})+\exp(\frac{-(u+v)^2}{4t})
\end{array}
\right)\otimes e^{-\mathbb A_t^2}\]
A similar formula, by exchanging the diagonal coefficients, 
gives the odd part of the heat kernel of $\bar{\mathsf  F}_t$ with 
boundary conditions \eqref{adboucon}. 
This expression and relation \eqref{asa} provide together an asymptotic  
explicite formula for $\bar K_+(t,u,v,y,z)$, the heat kernel of 
$\bar{\mathbb F}_t$. So we get the following formula, with respect to 
$C^\ell$-norm, for the supertrace 
density of $e^{-\bar{\mathbb F}_t}$ as a function of $t$ and $u$  
\begin{align*}
\bar K_+(t,u)&=\int_{N}\rm{str}\,\bar K_+(t,u,u,y,y)\,dy\\
&=-\frac{\rm{STr}\,e^{-\mathbb A_t^2}}{\sqrt{\pi t}}
e^{-\frac{u^2}{t}}+o(1)\in\Omega^*(B).
\end{align*}
Notice that in above formula, the integration is performed in 
linear space $\Omega^*(B)$ so it makes sens.
The case $\epsilon=-$ produces the same expression with the opposite sign. 
So by intergrating on $[0,1/2]$, with respect to $u$, we get the following 
asymptotic equality at $t=0$
\begin{equation*}
\bar K_{\epsilon_i}(t)=\int_0^{1/2}\bar K_{\epsilon_i}(t,u)\,du
\sim-\epsilon_i\frac{1}{2}\rm{STr}
\,e^{-\mathbb A_t^2}+o(1)\in\Omega^*(B).
\end{equation*}
In other word, the following relation holds in $\Omega^*(B)$ with 
$C^\ell$-topology
\begin{equation}\label{asd}
K_{\epsilon_i}(t)=-\epsilon_i\Ch(\mathbb A_{it})+o(1).
\end{equation}
\subsection{ Index density of APS conditions on family of half 
cylinders} 
Now consider the family of half cylinder $\R^{\geq0}\times N\to B$ with APS 
boundary condition. 
{\it We assume that $\dim\ker(A_b)$ does not depend on $b\in B$}.
Here $A_b$ indicates the Dirac operator on $N_b$. 
So vector spaces $\ker A_b$, put together, form the vector bundle 
$\ker A$ on $B$. 
We denote the lifting of this bundle, via $\pi$, by $\mathcal E_0$ 
which is again a vector bundle over $B$. 
In the same way and with notation of proposition \ref{ahah},
infinite dimensional vector spaces 
$E^\pm(b)$ together form an infinite dimensional bundles over $B$.
These bundles can be lifted to bundles $\mathcal E^\pm$ over $B$. 
These bundles are isomorphism via unitary operator $\mathcal U$. 
Let $\bar\nabla^\pm$ be connections on sections of 
$\mathcal E^\pm$ such that 
$\mathcal U\bar\nabla^\pm\mathcal U^{-1}=\bar\nabla^\mp$. 
These connections are assumed be constant along $\R^{\geq0}$, i.e. 
they are lifting of connections by $\pi$. 
We assume also a connection $\bar\nabla_0$ on $\mathcal E_0=\pi^*\ker A$. Put 
$\bar\nabla=\bar\nabla^-\oplus\bar\nabla_0\oplus\bar\nabla^+$ 
and consider the superconnection 
$\tilde{\mathbb B}
=\mathsf D+\bar\nabla\oplus\bar\nabla$. 
In view of relation \eqref{asvel}, heat operator of the rescaled 
supercurvature $\tilde{\mathbb F}_t$, up to a term of order $o(1)$, 
is equal to the heat operator of  
\begin{equation}
\mathsf  F_t:=\left(\begin{array}{cc}
-tD^*D&0\\0&-tDD^*
\end{array}
\right)+\left(\begin{array}{cc}
\bar\nabla^2&0\\0&\bar\nabla^2
\end{array}
\right)
\end{equation}
We shall prove that the supertrace of the heat operator $e^{-\mathsf  F_t}$ 
, when it is restricted to $\mathcal E^+\oplus\mathcal E^-$, 
vanishes when $t$ goes toward $0$. For this purpose, using the Voltera formula, 
we have 
\[(e^{-tD^*D+\bar\nabla^2})_{|\mathcal E^+}
=\sum_k(-1)^kI_k(tD^*D_{|\mathcal E^+},\bar\nabla_+^2)\]
where 
\[I_k(tD^*D_{|\mathcal E^+},\bar\nabla_+^2)
:=\int_{\triangle_k}e^{-s_0tD^*D}\bar\nabla_+^2e^{-s_1tD^*D}\bar\nabla_+^2
\dots e^{-s_{k-1}tD^*D}\bar\nabla_+^2e^{-s_ktD^*D}.\]
In other hand, from the first part of proposition \ref{ahah}, 
\[(e^{-tDD^*})_{|\mathcal  E^-}
=\mathcal U(e^{-tD^*D})_{\mathcal E^+}\mathcal U^{-1}.\]
This relation and $\bar\nabla_-=\mathcal U\bar\nabla_+\mathcal U^{-1}$ imply  
\[I_k(tD^*D_{|\mathcal E^-},\bar\nabla_-^2)
=\mathcal U I_k(tD^*D_{|\mathcal E^+},\bar\nabla_+^2)\mathcal U^{-1},\] so 
\[(e^{-tDD^*+\bar\nabla^2})_{|\mathcal E^-}
=\mathcal U (e^{-tD^*D+\bar\nabla^2})_{|\mathcal E^+}\mathcal U^{-1}.\]
A similar discussion, using again the proposition \ref{ahah}, implies 
the following relation in $\Omega^*(B)$ with $C^\ell$-topology
\[(e^{-tDD^*+\bar\nabla^2})_{|\mathcal E^+}
=\mathcal U (e^{-tD^*D+\bar\nabla^2})_{|\mathcal E^-}
\mathcal U^{-1}+o(1)\quad\text{ at }~t=0.\]
Therefore 
\[\rm{STr}\,
e^{-\bar{\mathbb F}_t}=\rm{Tr}\,(e^{-t\partial_u^2+\bar\nabla^2_0})_{|E_0} 
-\rm{Tr}\,(e^{-t\partial_u^2+\bar\nabla^2_0})_{|E_0}+o(1).\]
Since $\partial_u$ commutes with $\bar\nabla^2_0$, relation \eqref{laze} 
can be used to get
\[\rm{STr}_N\,e^{-\mathsf  F_{t}}(t,u)
=-\frac{e^{\frac{-u^2}{t}}}{\sqrt{\pi t}}
\rm{str}\,(e^{-\nabla_0^2})_{|\ker A}.\] 
Summarizing these discussion, we get the following relation in 
$\Omega^*(B)$ with $C^\ell$-topology
\begin{equation}\label{akhri}
\int_0^{\frac{1}{2}}\rm{STr}_N\,e^{-\tilde{\mathbb F}_t}(t,u)\,du
=-\frac{1}{2}\Ch(\ker A,\nabla_0)+o(1).
\end{equation}

\subsection{Index theorem for families of Dirac operators}
We recall that all differential operators and geometric structures 
we have on $F$ are of product form in collar neighborhood $U$, 
so they can be extended, smoothly, to double fibration $F\sqcup_{ F'}F$.
In opposite direction, let $\tilde\nabla$ be a connection on 
$\mathcal E\to F\sqcup_{ F'}F$ 
and let $\bar\nabla_i$ be the connections on half-cylinder fibrations 
discussed in previous subsections. Then 
\begin{equation}\label{conon}
\nabla:=\sum_if_1\bar\nabla_ig_1+f_2\tilde\nabla g_2,
\end{equation}
defines a connection on $\mathcal E\to B$. Here functions 
$f_1$, $f_2$, $g_1$ and $g_2$ are defined in subsection \ref{prforone}.
Let $\tilde{\mathbb F_t}$ be the rescaled supercurvature on double 
$F\sqcup_{ F'}F$ and let $\tilde K_t$ be the kernel of associated 
heat operator. We denote by $\bar K_{t\epsilon_i}$ and $K_{t\epsilon_i}$ the 
fundamental solutions of $e^{-\bar{\mathbb F}_t}$ and 
$e^{-\mathbb F_t}$ with respect to boundary condition $P^{\epsilon_i}$. 
At first we prove the following asymptotic formula at $t=0$
\begin{equation}\label{khastam}
K_{t\epsilon_i}=\sum_if_1\bar K_{t\epsilon_i}g_1+f_2\tilde K_tg_2+o(1).
\end{equation}
For this purpose we use relation \eqref{asvel} and Voltera formula to deduce 
\begin{gather}
e^{-\mathbb F_t}= e^{-t\mathsf D^2}+\sum_{k=1}^{\dim B}(-1)^k
I_k(t\mathsf D,\nabla\oplus\nabla)+o(1),\notag\\
I_k(t\mathsf D,\nabla\oplus\nabla):=\int_{\triangle_k}e^{-s_0t\mathsf D^2}
(\nabla\oplus\nabla)e^{-s_1t\mathsf D^2}(\nabla\oplus\nabla)\dots 
e^{-s_{k-1}t\mathsf D^2}(\nabla\oplus\nabla)e^{-s_kt\mathsf D^2}.\notag
\end{gather}
In this formula the boundary conditions are implicit in smooth family of 
heat operators $e^{-t\mathsf D^2}$. Following relation is the family 
version of \eqref{initst} which is clearly true
\[e^{-t\mathsf D^2}= f_1e^{-t\bar{\mathsf D}^2}g_1
+f_2e^{-t\tilde{\mathsf D}^2}g_2+o(1).\] 
In other hand   
$\nabla=\sum_if_1\bar{\nabla}^i g_1+f_2\tilde\nabla g_2$, 
since $f_1g_1=g_1$ and $f_2g_2=g_2$, we get 
\begin{align*}
I_k(t\mathsf D,\nabla\oplus\nabla)&=\sum_{\epsilon_i=\pm}f_1
I_k(\bar{\mathsf D}_{\epsilon_i},f_1(\bar{\nabla}^i\oplus\bar{\nabla}^i)g_1)g_1 
+f_2I_k(\tilde{\mathsf D},f_2(\tilde\nabla\oplus\tilde\nabla)g_2)g_2\\
&+\text{finite sum of operators of form }M(t)h
\,e^{-st\bar{\mathsf D}^2}\,kN(t)+o(1)
\end{align*}
In last line of above expression, $M(t)$ and $N(t)$ are vertical smoothing 
operators with differential form coefficients such that, far from diagonal, 
their kernels goes exponentially toward zero when $t$ goes to $0$, 
while $h$ and $k$ are smooth functions on $F$ such that 
$supp(hk)\subseteq[\frac{1}{4},1]\times F'\subset U$. So using relations 
\eqref{abed1} and \eqref{abed} we deduce 
$M(t)h\,e^{-st\bar{\mathsf D}^2}\,kN(t)=o(1)$ when $t\to0$. 
Using again the relation \eqref{abed1} and \eqref{abed} we have 
\[g_1e^{-t\bar{\mathsf D}^2}g_1=e^{-t\bar{\mathsf D}^2}g_1+o(1)~~\text{ and } 
~~f_1e^{-t\bar{\mathsf D}^2}f_1=f_1e^{-t\bar{\mathsf D}^2}+o(1).\]
Therefore  
\[f_1I_k(\bar{\mathsf D}_{\epsilon_i},
f_1(\bar{\nabla}^i\oplus\bar{\nabla}^i)g_1)g_1
=f_1I_k(\bar{\mathsf D}_{\epsilon_i},\bar{\nabla}^i\oplus\bar{\nabla}^i)
g_1+o(1)\] 
which complete the proof of relation \eqref{khastam} by considering 
\eqref{asvel}.
\newline Now we take the supertrace of relation \eqref{khastam}. 
Clearly the contribution of interior term is zero. Contribution of 
boundary terms, $\bar K_{t\epsilon_i}g_1$ on $[1/2,1]$ go to $0$ when 
$t$ goes toward $0$ (see relation \eqref{abed1} and \eqref{abed}). 
Since $f_1=g_1=1$ on $[0,1/2]$, using relations \eqref{asd} and \eqref{akhri} 
we obtain the following asymptotic formula at $t=0$ 
\begin{equation}\label{erfd}
\rm{STr}\,(e^{-\mathbb F_t^\epsilon})
=-\frac{1}{2}\sum_{\epsilon_i=\pm}\epsilon_i\,\Ch(\mathbb A_{it})
-\frac{1}{2}\sum_{\epsilon_i=0}\Ch(\ker A_i,\nabla_0)+o(1)\in\Omega^*(B)
\end{equation}
All involved differential forms in above expressions are closed, so we can 
regard this relation in $H^*_{dr}(B)$. 
According to proposition \ref{pro1} and its analogue for 
superconnection $\mathbb A$, the class of Chern forms do not depend 
on parameter $t$. So the image of term $o(1)$ in $H^*_{dr}(B)$ vanishes and we 
obtain the following equality 
\[\Ch[\ind(D,P^\epsilon)]=-\frac{1}{2}\sum_{\epsilon_i=+}\Ch[\ind A_i] 
+\frac{1}{2}\sum_{\epsilon=-}\Ch[\ind A_i]
-\frac{1}{2}\sum_{\epsilon_i=0}\Ch(\ker A_i)\in H^*_{dr}(B).\]
This complete the proof of the theorem \ref{thm2}.
\begin{rem}
Above proof is based on assumption that $\dim\ker A_b$ 
is independent of $b\in B$. This assumption is satisfied for some interesting 
cases, eg. for family of signature operators or, when fibers have a metric with 
positive scalar curvature, for standard Dirac operator twisted by flat 
vector bundles. However this assumption may be removed by considering 
smooth perturbations of boundary operators, or putting more general 
spectral boundary condition by means of spectral projections introduced 
in \cite{MePi}. In the later case the standard tool for 
analysing the heat kernel will be the Melrose's b-calculus for family. 
\end{rem}


\providecommand{\bysame}{\leavevmode\hbox to3em{\hrulefill}\thinspace}
\providecommand{\MR}{\relax\ifhmode\unskip\space\fi MR }
\providecommand{\MRhref}[2]{%
  \href{http://www.ams.org/mathscinet-getitem?mr=#1}{#2}
}
\providecommand{\href}[2]{#2}

\end{document}